\documentclass[12]{article}
\usepackage[dvips]{graphicx}
\def\no{\noindent}

\def\bs{{\bf s}}

\def\bA{{\bf A}}
\def\bB{{\bf B}}

\def\bx{{\bf x}}

\def\bn{{\bf n}}
\def\ba{{\bf a}}
\def\bb{{\bf b}}

\def\be{{\bf e}}

\def\C{{C\kern-.647em I}}

\def\R{I\!\!R}

\def\beq{\begin{equation}}
\def\eeq{\end{equation}}

\def\w{\wedge}

\def\no{\noindent}

\def\bB{{\bf B}}

\def\C{{I\! \! \! C}}

\def\no{\noindent}

\def\no{\noindent}

\def\G{I\!\!\!G}

\def\d{\partial}
\def \xx {{\bf x}}

\def \G {{\cal G}}

\begin{document}
   
\title{Fundamental Theorem of Calculus}

\author{Garret Sobczyk   \\ 
Universidad de Las Am\'ericas - Puebla,\\ 72820 Cholula, Mexico, \\
Omar Sanchez\\ University of Waterloo, \\ Ontario, N2L 3G1 Canada}

\maketitle

\begin{abstract} A simple but rigorous proof of the {\it Fundamental Theorem of Calculus} is given
in geometric calculus, after the basis for this theory in geometric algebra has been explained. 
Various classical examples of this theorem, such as the Green's and Stokes' theorem are discussed,
as well as the new theory of monogenic functions, which generalizes the concept of an analytic function
of a complex variable to higher dimensions.
  
 
  
 \end{abstract}

\section*{Introduction}

  Every student of calculus knows that integration is the inverse operation to differentiation.
Things get more complicated in going to the higher dimensional analogs, such as {\it Green's theorem} in the
plane, {\it Stoke's theorem}, and {\it Gauss' theorem} in vector analysis. But the end of the matter
is nowhere in sight when we consider differential forms, tensor calculus, and differential geometry,
where important analogous theorems exist in much more general settings. The purpose of this article is
to express the fundamental theorem of calculus in geometric algebra, showing that 
the multitude of different forms of this theorem fall out of a single theorem,
and to give a simple and rigorous proof. An intuitive proof of the fundamental theorem
in geometric algebra was first given in \cite{H68}. The basis of this theorem was further explored in
 \cite{Simsob}, \cite{omar}. We believe that the geometric algebra version of this
important theorem represents a significant improvement and generalization over other forms found in the literature,
such as for example in \cite{flanders}.

Geometric algebra was introduced by William Kingdon Clifford in 1878 as a generalization and unification
of the ideas of Hermann Grassmann and William Hamilton, who came before him 
\cite{Clifford}, \cite{Crowe1985}. Whereas Hamilton's quaternions are
well known and Grassmann algebras are known in various guises such as the exterior algebras of differential 
forms \cite{RS}, and anti-symmetric tensors, the full 
generality and utility of Clifford's geometric algebras are just beginning to be
appreciated \cite{EBGS}, \cite{DDL}. We hope that this paper will attract the attention of the 
scientific community to the generality and beauty of this formalism. 
What follows is a short introduction to the basic ideas of geometric algebra
that are used in this paper. 

\section*{Geometric algebra}

   Let $\{\be_i\}_{i=1}^n$ be the standard orthonormal basis of the $n$-diminsional euclidean space
 $\R^n$. Vectors $\ba, \bb \in \R^n$ can be represented in the alternative forms
   \[ \ba=(a^1 \ a^2 \ \cdots \ a^n) = \sum_i a^i \be_i \ \ {\rm and} \ \
    \bb=(b^1 \ b^2 \ \cdots \ b^n) = \sum_i b^i \be_i,\]
  and their {\it inner product} is given by 
    \[\ba \cdot \bb=|\ba||\bb|\cos \theta = \sum_i a^i b^i  \]
 where $\theta$ is the angle between them.     

The history of mathematics begins with the invention of the natural numbers \cite{D1954}, and
continues with successive extensions of the number concept \cite{St}.      
The {\it associative geometric algebra} $\G_n=\G_n(\R^n)$ of the 
euclidean space $\R^n$ can be thought of as the {\it geometric
extension} of the real number system $\R$ to include $n$ new {\it anticommuting} square roots of unity, which we
identify with the basis vectors $\be_i$ of $\R^n$, and their geometric products. 
This is equivalent to saying that 
each vector $\ba \in \R^n$ has the property that 
  \[  \ba^2 = \ba \ba = \ba\cdot \ba = |\ba|^2, \]
where $|\ba|$ is the euclidean length or {\it magnitude} of the vector $\ba$.  

The fundamental {\it geometric product} $\ba \bb$ of the vectors $\ba, \bb \in \R^n \subset \G_n$ can be decomposed
into the {\it symmetric} inner product, and the {\it anti-symmetic} outer product
  \beq \ba \bb = \frac{1}{2}(\ba \bb + \bb \ba) + \frac{1}{2}(\ba \bb - \bb \ba)=\ba \cdot \bb + \ba \w \bb
  =|\ba||\bb|e^{i \theta}  ,    \label{vecgp} \eeq
where $\ba \cdot \bb =  \frac{1}{2}(\ba \bb + \bb \ba)$ is the inner product and 
\[ \ba \w \bb =  \frac{1}{2}(\ba \bb - \bb \ba) = |\ba||\bb|i \sin \theta \]
 is the {\it outer product} of the
vectors $\ba$ and $\bb$. The outer product $\ba \w \bb$
is given the geometric interpretation of a {\it directed plane segment} or {\it bivector} and characterizes the
direction of the plane of the subspace of $\R^n$ spanned by the vectors $\ba$ and $\bb$. 
The {\it unit bivector} $i$ {\it orients} the plane of the vectors $\ba$ and $\bb$, and has
the property that $i^2=-1$. See Figure 1. 
The last equality in (\ref{vecgp}) expresses the beautiful Euler formula for the geometric product
of the two vectors $\ba$ and $\bb$.
\begin{figure}
\begin{center}
\includegraphics[scale=.80]{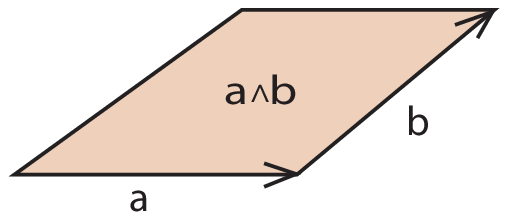}
\includegraphics[scale=.70]{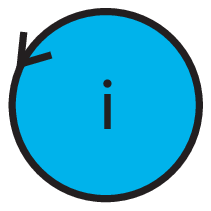}
\caption{The bivector $\ba\w \bb$ is obtained by sweeping 
the vector $\bb$ out along the vector $\ba$. Also shown is the unit bivector $i$ in the plane of
$\ba$ and $\bb$.}
\end{center}
\end{figure} 

The real $2^n$-dimensional geometric algebra $\G_n=\G_n(\R^n)$ has the {\it standard basis}
  \[  \G_n={\rm span}\{ 1,\be_{\lambda_1},\be_{\lambda_2}, \dots , \be_{\lambda_n}\}, \]
where the $\pmatrix{n \cr k}$ {\it $k$-vector} basis elements of the form $\be_{\lambda_k}$ 
are defined by  
  \[ \be_{\lambda_k}=\be_{i_1\ldots i_k}=\be_{i_1} \cdots \be_{i_k}  \]    
for each $\lambda_k=i_1 \cdots i_k$ where $1 \le i_1 < \cdots < i_{k}\le n$. 
Thus, for example,
the $2^3$-dimensional geometric algebra of $\R^3$ has the standard basis
  \[  \G_3={\rm span}\{ 1,\be_1,\be_2,\be_3,\be_{12},\be_{13},\be_{23},\be_{123}\}. \]
We denote the {\it pseudoscalar} of $\R^n$ by the special symbol $I=\be_{1\ldots n}$. The
pseudoscalar gives a unique orientation to $\R^n$ and $I^{-1}=\be_{n \ldots 1}$.

   Let $\ba$ be a vector, and $\bB_k$ be a $k$-vector in $\G_n$ where $k \ge 1$. Just as for (\ref{vecgp}),
we decompose the product $\ba \bB_k$ into symmetric and anti-symmetric parts,
   \beq  \ba \bB_k =  \frac{1}{2}(\ba \bB_k + (-1)^{k+1}\bB_k \ba) +
    \frac{1}{2}(\ba \bB_k - (-1)^{k+1} \bB_k \ba)=\ba \cdot \bB_k + \ba \w \bB_k ,
                       \label{kvecgp} \eeq 
where 
  \[\ba \cdot \bB_k = \frac{1}{2}(\ba \bB_k + (-1)^{k+1}\bB_k \ba)=<\ba \bB_k>_{k-1}\]
   is a $(k-1)$-vector and
 \[ \ba \w \bB_k =  \frac{1}{2}(\ba \bB_k - (-1)^{k+1} \bB_k \ba)=<\ba \bB_k>_{k+1}\]
  is a $(k+1)$-vector. More generally, if $\bA_r$ and $\bB_s$ are $r$- and $s$-vectors of $\G_n$, 
  where $r>0, s>0$
  we define $\bA_r \cdot \bB_s=<\bA_r \bB_s>_{|r-s|}$ and $\bA_r \w \bB_s=<\bA_r \bB_s>_{r+s}$. 
  In the exceptional cases when $r=0$ or $s=0$, we define $\bA_r\cdot \bB_s=0$ and 
  $\bA_r\w \bB_s =\bA_r \bB_s$.
  
  One very useful identity
     \[ \ba\cdot (\bA_r\w \bB_s)=(\ba \cdot \bA_r)\w B_s+(-1)^r A_r \w (\ba \cdot \bB_s)=
       (-1)^{r+s+1}(\bA_r\w \bB_s)\cdot \ba     \]
gives the {\it distributive law} for the inner product of a vector over the outer product of an
$r$- and an $s$-vector. For an $n$-vector $\bA_n\in \G_n$, the {\it determinant function} is
defined by $\bA_n=\det(\bA_n)I$ or $\det(\bA_n)=\bA_n I^{-1}$.  
   
      We have given here little more than a list of some of the basic algebraic 
identities that will use in this paper. However, the subject has become highly developed in recent years. 
We suggest the following references \cite{H/S}, \cite{LP97}.
In addition, the following links provide good on-line references \cite{links}.
 
  \section*{Calculus on a $k$-surface in $\cal M$}

Let $R=\times_{i=1}^k[a^i,b^i]$ be a $k$-rectangle in $\R^k$. We denote the points $\bs \in R$ by
  the {\it position vectors} $\bs =\sum_i s^i \be_i$ where $\{\be_i\}_{i=1}^k$ 
  is the standard orthonormal basis of $\R^k$
  and $a^i\le s^i \le b^i$.
The {\it boundary} of the rectangle is the {\it $k$-chain}
\begin{displaymath}
\beta(R) = \oplus_{i=1}^k (R^i_+ \oplus R^i_-)
\end{displaymath}
where the faces $R^i_\pm$ are defined by $R^i_+ =R(s^i=b^i)$ and $R^i_-=R(s^i=a^i)$, respectively. 

Let $\cal M$ be a $m$-surface in $\R^n$ which is of class $C^p$, where $p\ge 2$ and $m\le n$.
Such a surface can be defined in the usual way using the
{\it charts} and {\it atlases} of differential geometry, but we are not interested in such details here.
Instead, we are interested in studying the properties of $m$-dimensional {\it rectangular patches} on
the surface $\cal M$, and then the surfaces obtained by piecing together such patches, in much the
same way that rectangular strips are used to define the {\it Riemann integral} for the
area under a curve in elementary
calculus.

\bigskip

\no {\bf Definition}
A set $M\subset {\cal M}\subset \R^n$ is a {\it rectangular $k$-patch} of class $C^p$ ($p\geq 1$) in $\cal M$ 
if $M=\xx(R)$, where $\xx\colon R\to  M$ is a proper ($\xx^{-1}:M \to R$ is continuous),
 regular $k$-patch of class $C^p$. Of course, $1\le k\le m$.

\bigskip

As an example, we give the graph of the image of the $2$-rectangle $R=[0,1]\times [0,1]$ where 
\[ \bx(s^1,s^2)=(s^1,s^2,\frac{1- \sin((s^1)^2)}{2}-3s^2), \] 
see Figure 2.

\begin{figure}
\begin{center}
\includegraphics[scale=1.40]{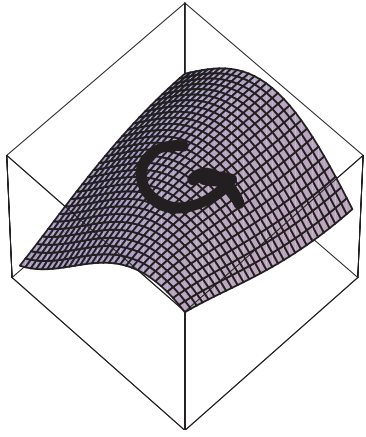}%{fig2a.eps}
\caption{The oriented $2$-patch $M=\bx(R)$ is the image of the $2$-square $R=[0,1]\times [0,1]$. The
orientation is shown by the curved arrow.}
\end{center}
\end{figure} 

The {\it boundary} of $M$ is the $k$-chain 
\begin{displaymath}
\beta(M) =\xx(\beta(R))= \oplus_{i=1}^k (M^i_+ \oplus M^i_-)
\end{displaymath}
where the {\it faces} $M^i_\pm=\xx(R^i_\pm)$, respectively, for $i=1,\dots,k$. 
The tangent vectors $\{\xx_i\}_{i=1}^k$ at a point $\bx\in M$, defined by
\begin{equation}
\xx_i=\frac{\d \xx}{\d s^i}, \qquad i=1,\dots,k,
\end{equation}
make up a local basis of the {\it tangent space} at the point $\bx\in M$, 
and for a proper, regular patch they are linearly independent and continuous on $M$.

The tangent $k$-vector $\xx_{(k)}$ at a point $\bx\in M$ is the  defined by 
\begin{equation}
\xx_{(k)}=\xx_1\wedge \cdots \wedge \xx_k,
\end{equation}
and for a proper, regular patch is continuous and doesn't vanish at any point $\bx \in M$.
The {\it reciprocal basis} $\{\bx^i\}_{i=1}^k$ to $\{\xx_i\}_{i=1}^k$ at the point $\bx \in M$
is defined by
  \[ \bx^i =(-1)^{i-1}(\w_{j=1,j\ne i}^k \bx_i)\cdot \bx_{(k)}^{-1} \]
and satisfies the relations $\bx^i\cdot \bx_j = \delta^i_{\ j}$. For example for $i=1$,
  \[ \bx^1 = (\bx_2 \w \cdots \w \bx_k)\cdot  \bx_{(k)}^{-1} . \]  

The oriented tangent $(k-1)$-vector on each face $M^i_\pm$ of the boundary of $M$ is
defined by 
\begin{equation}
\xx_{(k-1)}=\pm \xx_{(k)}\xx^i,
\end{equation}
respectively.  The vector $\xx^i$ is reciprocal to $\xx_i$, and $\pm \xx^i$ defines
the directions of the {\it outward normals} to the faces $M^i_\pm$ at the points $\bx\in M^i_\pm$, respectively.

Let $f,g\colon M\to \G_n$ be continuous $\G_n$-valued 
functions on the patch $M \subset {\cal M} \subset \R^n$. More precisely, we should define $f$ and $g$ to
be continuous on an open set ${\cal S}$ such that $M \subset {\cal S} \subset {\cal M}$, since we will
later want $f$ and $g$ to be defined on an open set containing all of $\cal M$.

\bigskip

\no {\bf Definition:} The directed integral over $M$ is defined by
\begin{equation}
\int_{M}gd\xx_{(k)}f=\int_{R}g\xx_{(k)}fds^{(k)} ,
\end{equation}
where $ds^{(k)}=ds^1\cdots ds^k$.
We also write $\frac{d \bx_{(k)}}{ds^{(k)}}=\frac{\partial \bx}{\partial s^1}\w \cdots 
\w \frac{\partial \bx}{\partial s^k}$.

\bigskip

\no {\bf Definition:} The directed integral over the boundary is
\begin{displaymath}
\int_{\beta(M)}gd\xx_{(k-1)}f=
\sum_{i=1}^n \int_{M^i_+ \oplus M^i_-}gd\xx_{(k-1)}f
\end{displaymath}
\begin{displaymath}
=\sum_{i=1}^n  \int_{R^i_+ \oplus R_i^-}g\xx_{(k-1)}f ds^{(k-1)}
\end{displaymath}
where $ds^{(k-1)}=\frac{ds^{(k)}}{ds^i}=ds^1\cdots\widehat ds^i \cdots ds^k$ in $R^i_\pm$.

\bigskip

The {\it directed content} of $M$ is defined by
\begin{displaymath}
{\cal D}(M)=\int_{M}d\xx_{(k)}
\end{displaymath}
A direct consequence of the fundamental theorem will be that
\begin{displaymath}
{\cal D}(\beta(M))=0.
\end{displaymath}

The {\it two sided} vector derivative $g(\bx)\d_\bx f(\bx)$  of $f(\bx)$ and $g(\bx)$ on $M=\bx(r)$
 at the point $\bx\in M$ is defined indirectly by
\begin{equation}
g(\bx)\d_\bx f(\bx) 
=\sum_{i=1}^k(\xx_i \cdot \d_\bx)(\dot g \xx^i \dot f)=\sum_{i=1}^k\frac{\partial}{\partial s^i}(\dot g \xx^i \dot f),
\end{equation}
where the deriviatives act only on the dotted arguments.
 Note in the above definition that we are using the {\it chain rule} 
  \[  \partial = \sum_{i=1}^k \xx^i \frac{\partial \bx}{\partial s^i}\cdot
    \partial_\bx = \sum_{i=1}^k \xx^i \frac{\partial}{\partial s^i}
   . \]
It follows that the reciprocal vectors $\xx^i$ can be expressed in
terms of the vector derivative by $\xx^i=\partial s^i$ for
$i=1, \cdots, k$. The {\it Libnitz product rule} for partial derivatives gives 
\begin{displaymath}
g\d f=\dot g \d f+ g \d \dot f.
\end{displaymath}

\bigskip

\no {\bf Lemma:} Let $M$ be a rectangular $k$-patch of class $C^2$ in $\R^n$. Then 
   $$\sum_{i=1}^j \frac{\partial}{\partial s^i} \xx_{(j)}\cdot \bx^i=0,$$
 for $j=1,\ldots, k$. In the case when $j=k$, the $` ` \cdot "$ can be removed. 
 
 {\bf Proof:} For $j=1$, clearly  $\frac{\partial}{\partial s^1} \bx_1 \cdot \bx^1=
 \frac{\partial}{\partial s^1} 1=0$. 
 We now inductively assume that the lemma is true for all $j<k$, and calculate
   \[ \sum_{i=1}^k \frac{\partial}{\partial s^i} \xx_{(k)}\cdot \bx^i = 
    \sum_{i=1}^k \frac{\partial}{\partial s^i} [(\xx_{(k-1)}\w \bx_k)\cdot \bx^i] \]
    \[ = \frac{\partial}{\partial s^k}\bx_{(k-1)}-
    \sum_{i=1}^{k-1}\frac{\partial}{\partial s^i}[(\bx_{(k-1)}\cdot \bx^i)\w \bx_k] \] 
    \[ = \bx_{(k-1),k} -  \sum_{i=1}^{k-1}(\bx_{(k-1)}\cdot \bx^i)\w \bx_{k,i} 
               =0. \] 
  In the last step, we have used the fact that partial derivatives commute so that 
  \[ \bx_{k,i}=\frac{\partial \bx_k}{\partial s^i} = \frac{\partial \bx_i}{\partial s^k} = \bx_{i,k} . \] 
  
      \hfill $ \framebox{} $             

\bigskip

\no {\bf Theorem}[Fundamental Theorem] Let 
$M$ be a rectangular $k$-patch of class $C^2$, and $f,g\colon M\to \G_n$ of class $C^1$. Then
\begin{equation}
\int_{M}g d\xx_{(k)}\d f=\int_{\beta(M)}gd\xx_{(k-1)}f
\end{equation}

\bigskip

\no {\bf Proof:}
\begin{displaymath}
\int_{\beta(M)}gd\xx_{(k-1)}f=
\end{displaymath}
\begin{displaymath}
=\sum_{i=1}^k \int_{M^i_+\oplus M^i_-}gd\xx_{(k-1)}f=\sum_{i=1}^n  \int_{R^i_+ \oplus R^i_-}g\xx_{(k-1)}f ds^{(k-1)}
\end{displaymath}
\begin{displaymath}=\sum_{i=1}^n  \int_{R^i_+\oplus R_i^-}\Big[
\int_{a_i}^{b_i}\frac{\partial}{\partial s^i}(g\xx_{(k)}\bx^i f) ds^i\Big]\frac{ds^{(k)}}{ds^i} \]
\[=\sum_{i=1}^n  \int_{R}\Big[ \frac{\partial}{\partial s^i}(\dot g \xx_{(k)} \bx^i \dot f) \Big] ds^{(k)}
\end{displaymath}
\begin{displaymath}
=\sum_{i=1}^k \int_{M}(\xx_i\cdot \d) \dot g {d\xx}_{(k)}  {\xx^i} \dot f=\int_{M} gd\xx_{(k)}\d f.
\end{displaymath}
 \hfill $ \framebox{} $        

\bigskip

   Choosing $g=f=1$, the constant functions, in the Fundamental Theorem immediately gives
   
   \bigskip
   
   \no {\bf Corollary} Let $M$ be a rectangular $k$-patch of class $C^2$, then
     \[ {\cal D}(M)=\int_{\beta(M)} d\bx_{(k-1)}=0. \] 
     
     \bigskip

The fundamental theorem can now be easily extended in the usual way 
to any simple $k$-surface ${\cal S}$ in $\cal M$ 
obtained by glueing together $k$-rectangular blocks, making sure that proper orientations are respected on
their boundaries.

% Figure 3 shows how this is done for $2$-dimensional blocks when $k=2$. 

\section*{Classical theorems of integration}

Let $\bx(s)$, for $s\in M=[a,b]$ be a $1$-dimensional curve on the $C^2$ surface $\cal M$, and let
$f,g$ be $C^1$ functions mapping $M \to \G_n$. Then the fundamental theorem gives
 \[  \int_{M}g d\xx \cdot \d f=
  \int_{\beta(M)}gd\xx_{(0)}f = g(\bb)f(\bb)- g(\ba)f(\ba), \]
 where $\ba=\bx(a)$ and $\bb=\bx(b)$. Of course, this gives the condition of when the integral along a curve
 connecting the points $\ba, \bb \in {\cal M}$ is independent of the path. Note that the
 vector derivative differentiates both to the left and to the right.

For a $2$-dimensional simple surface ${\cal S}$ embedded in a higher dimensional $k$-manifold  
${\cal M} $, with boundary $\beta({\cal S})$, the fundamental theorem gives
  \beq  \int_{\cal S}g d\xx_{(2)}\cdot \d f=\int_{\beta({\cal S})}gd \xx f . \label{funthm2} \eeq
This is a general integration theorem that includes 
the classical Green's theorem in the plane and Stokes' theorem, among many others, as special cases \cite{outline}.
 If ${\cal M}=\R^2$, then the $2$-surface ${\cal S}$ lies in the plane, and the fundamental
 theorem (\ref{funthm2}) takes the form
  \[   \int_{\cal S}g d\xx_{(2)} \d f=\int_{\beta({\cal S})}gd \xx f .  \]
  
   Choosing
$g(\bx)=1$ and $f=f(\bx)$ to be a $C^1$ vector-valued function, and identifying scalar
and bivector parts of the last equation, gives
  \[   \int_{\cal S} d\xx_{(2)} \d \w f=\int_{\beta({\cal S})} d \xx \cdot f \ \ {\rm and} \ \  
  \int_{\cal S} d\xx_{(2)} \d \cdot f=\int_{\beta({\cal S})} d \xx \w f, \]
 or, equivalently
  \[   \int_{\cal S} d\xx_{(2)} \d \w f=\int_{\beta({\cal S})} d \xx \cdot f \ \ {\rm and} \ \  
  \int_{\cal S} d\xx_{(2)} \d \cdot f=\int_{\beta({\cal S})} d \xx \w f, \]
  which are two quite different forms of the standard Green's Theorem \cite[p. 110]{outline}.  
  
   If ${\cal S}$ is a simple $2$-surface in
  ${\cal M}=\R^3$ having the simple closed curve $\beta({\cal S})$ as boundary, then 
  taking the scalar parts of (\ref{funthm2}), and utilizing the {\it cross product} of
  standard vector analysis, gives
    \[ \int_{\cal S}( \d \times f)\cdot \hat \bn |d\xx_{(2)}| =\int_{\beta({\cal S})} f \cdot d \xx   \]
  which is Stokes' Theorem.   
   $g=1$ and $f=f(\bx)$ is a vector-valued function, then the
   above integral becomes equivalent to the Stokes' Theorem for a simple surface in $\R^3$.   

 If ${\cal S}$ is a simple closed $2$-surface in
  ${\cal M}=\R^3$, $g=1$ and $f=f(\bx)$ is a vector-valued function, then the
  fundamental theorem gives  
   \[  \int_{\cal S} d\xx_{(3)} \d \cdot f=\int_{\beta({\cal S})}d \xx_{(2)} \w f, \]
 which is equivalent to the Gauss' Divergence Theorem. This becomes more evident when we multiply
 both sides of this integral equation by $I^{-1}=\be_{321}$, giving
    \[  \int_{\cal S} \d \cdot f \ |d\xx_{(3)}|=\int_{\beta({\cal S})}  \hat \bn \cdot f\ |  d \xx_{(2)}|, \]
where $\hat \bn$ is the unit outward normal to the surface element $d \xx_{(2)}$ at the point $\bx \in
 \beta({\cal S})$. There are many other forms of the classical integration theorems, and all of them
 follow directly from our fundamental theorem.  
 
Over the last several decades a powerful {\it geometric analysis} has been developed in
geometric algebra. One of the most striking features of this new theory is the generalization of
the concept of an analytic function to that of a {\it monogenic function} \cite{ca}. 

\bigskip

\no {\bf Definition:} A geometric-valued $C^1$ function $f(\bx)$ 
is said to be {\it monogenic} on ${\cal S}\subset {\cal M}$ if $\d_\bx f=0$ for
all $\bx \in {\cal S}$.

\bigskip

Let ${\cal M}=\R^n$ and let $\cal S$ be the $(n-1)$-dimensional simple closed boundary containing
an open region $M$ of $\cal M$. Let $f=f(\bx)$ be a $C^1$ geometric-valued function defined on $M$
and its boundary ${\cal S}$. Then for any point $\bx^\prime \in M$, 
  \[  f(\bx^\prime) = \frac{1}{\Omega} \int_M |d\bx_{(n)}|\frac{(\bx^\prime-\bx)}{|\bx^\prime
  -\bx|^n}\partial_{\bx}f(\bx)
          -\frac{1}{\Omega} \int_{\cal S} |d\bx_{(n-1)}|\frac{(\bx^\prime-\bx)}{|\bx^\prime-\bx|^n}\hat \bn f(\bx), \]
 where $\hat \bn$ is the unit outward normal to the surface ${\cal S}$ and $\Omega = \frac{2 \pi^{n/2}}
 {\Gamma(n/2)}$ is the area of the $(n-1)$-dimensional unit sphere in $\R^n$. This important result
 is an easy consequence of the fundamental theorem by choosing $g(\bx)$ to be the {\it Cauchy kernel function}
 $g(\bx^\prime)=  \frac{(\bx^\prime-\bx)}{|\bx^\prime -\bx|^n} $, 
 which has the property that $\partial_{\bx^\prime}g(\bx^\prime)=0$ for all $\bx^\prime \ne \bx \in \R^n$.
 For a monogenic function $f(\bx)$, we see that     
   \[  f(\bx^\prime) = -\frac{1}{\Omega} \int_{\cal S} |d\bx_{(n-1)}|\frac{(\bx^\prime-\bx)}{|\bx^\prime-\bx|^n}\hat \bn f(\bx). \]
   This means that a monogenic function $f(\bx)$ in a region $M\subset \R^n$ 
   is completely determined by its values on its boundary $\cal S$. This, of course, generalizes
   Cauchy's famous integral theorem for analytic functions in the complex number plane \cite{A79}.
   
 \section*{The vector derivative in $\R^n$}
 
   General formulas for the vector derivative on a $k$-surface would take us deep into differential geometry
   \cite{H/S}. Here, we will only give a few formulas for the vector derivative on $\R^n$, which were utilized in
   the derivation of the generalized Cauchy integral formula given in the last section. The vector
   derivative $\partial_\bx$ on $\R^n$ can be most simply expressed in terms of the standard basis,
   $\partial=\partial_\bx = \sum_{i=1}^n \be_i \frac{\partial}{\partial x^i}$. Where the position vector $\bx\in \R^n$
   has the coordinate form 
     \[ \bx=(x^1 \ \cdots \  x^n )= \sum_{i=1}^n x^i \be_i . \] 
   
     The following formulas are easily verified:
     
     \bigskip
     
     \no 1. $ \partial \cdot \bx = n$ and $\partial \w \bx =0$, as easily follows from 
     $\partial \ \bx = \partial \cdot \bx + \partial \w \bx = n$,          
     
     \bigskip
     
     \no 2. $\ba \cdot \partial  \ \bx = \ba = \partial \ \bx\cdot \ba$, as easily follows from 
        $\ba\cdot (\partial \w \bx) =\ba\cdot \partial \ \bx - \partial \ \bx\cdot \ba =0, $
      
        \bigskip  
               
     \no 3. $\partial \ \bx^2 = \partial \ \dot \bx\cdot \bx +\partial \ \bx\cdot \dot \bx = 2\bx$,
     
     \bigskip
     
     \no 4. $\partial |\bx|= \partial (\bx^2)^{1/2} =\frac{1}{2}(\bx^2)^{-1/2}2 \bx = \hat \bx$,
     
     \bigskip
     
     \no 5. $\partial |\bx|^k = k|\bx|^{k-2}\bx$,
     
     \bigskip
     
     \no 6. $\partial \frac{\bx}{|\bx|^k}=\frac{n-k}{|\bx|^k}$,
     
     \bigskip
     
     \no 7. $\partial \log |\bx|= \frac{\bx}{|\bx|^2}=\bx^{-1}$.

 \section*{Acknowledgements} 
 
     The first author thanks Dr. Guillermo Romero, Academic Vice-Rector, and Dr. Reyla Navarro, Chairwomen of
   the Department of Mathematics, at the Universidad de Las Americas for continuing 
   support for this research. He and is a member of SNI 14587.

\end{document}